\documentclass[a4paper,11pt,reqno]{amsart}
\usepackage[T1]{fontenc}
\usepackage[utf8]{inputenc}
\usepackage{lmodern}
\usepackage{amsmath, amsthm, amssymb,amscd, mathrsfs, amsfonts, mathtools}
\usepackage{hyperref}
\hypersetup{backref=true}
\usepackage{euler}

\usepackage{afterpage}

\usepackage[all]{xy}
\usepackage{todonotes}
\usepackage{xcolor}

\usepackage[lite]{amsrefs}

\renewcommand{\PrintDOI}[1]{\href{http://dx.doi.org/\detokenize{#1}}{doi: \detokenize{#1}}%
  \IfEmptyBibField{pages}{, (to appear in print)}{}}

\def\commutatif{\ar@{}[rd]|{\circlearrowleft}}

\newcommand{\eq}[1][r]
   {\ar@<-3pt>@{-}[#1]
    \ar@<-1pt>@{}[#1]|<{}="gauche"
    \ar@<+0pt>@{}[#1]|-{}="milieu"
    \ar@<+1pt>@{}[#1]|>{}="droite"
    \ar@/^2pt/@{-}"gauche";"milieu"
    \ar@/_2pt/@{-}"milieu";"droite"}

\def\dar[#1]{\ar@<2pt>[#1]\ar@<-2pt>[#1]}
  \entrymodifiers={!!<0pt,0.7ex>+} 

\newcommand{\bigon}[4][r]{
    \ar@/^1pc/[#1]^{#2}_*=<0.3pt>{}="HAUT"
    \ar@/_1pc/[#1]_{#3}^*=<0.3pt>{}="BAS"
    \ar@{=>} "HAUT";"BAS" ^{#4}
  }

\newcommand{\bigons}[6][r]{  
    \ar@/^2pc/[#1]^{#2}_*=<0.3pt>{}="HAUT"
    \ar@{}    [#1]     ^*=<0.3pt>{}="MILIEUHAUT"
                       _*=<0.3pt>{}="MILIEUBAS"
    \ar[#1]_(0.3){#3}                  
    \ar@/_2pc/[#1]_{#4}^*=<0.3pt>{}="BAS"
    \ar@{=>} "HAUT";"MILIEUHAUT" ^{#5}
    \ar@{=>} "MILIEUBAS";"BAS" ^{#6}
  }

\newtheorem{thm}{Theorem}[section]
\newtheorem{pro}[thm]{Proposition}
\newtheorem{lem}[thm]{Lemma}

\theoremstyle{definition}
\newtheorem{df}[thm]{Definition}
\newtheorem{dfpro}[thm]{Definition and Proposition}

\theoremstyle{remark}
\newtheorem{rmk}[thm]{Remark}

\newtheorem{ex}[thm]{Example}

\allowdisplaybreaks

\newcommand{\Z}{\mathbb{Z}}

\newcommand\rTo{\longrightarrow}

\newcommand\mto{\longmapsto}

\newcommand\rRack{\triangleright}
\newcommand\Rack{\diamond}

\let\fr\mathfrak
\let\cal\mathcal

\let\bb\mathbb

\def\s{\sigma}

\hypersetup{
  colorlinks = true,
  urlcolor = blue,
  linkcolor = blue,
  citecolor = red,
  pdfauthor = {Elhamdadi, M.,},
  pdfkeywords = {An Algebraic Approach to Quandles: Recent development },
  pdftitle = {Elhamdadi, M. - An Algebraic Approach to Quandles: A Survey},
  pdfsubject = {racks, quandles},
  pdfpagemode = UseNone
}

\title{ A Survey of Racks and Quandles: Some recent developments. }

\author{Mohamed Elhamdadi} 
\address{Department of Mathematics, 
University of South Florida, Tampa, FL 33620, U.S.A.} 
\email{emohamed@math.usf.edu} 

\begin{document}

\maketitle

\begin{abstract}
This short survey contains some recent developments of the algebraic theory of racks and quandles.  We report on some elements of representation theory of quandles and ring theoretic approach to quandles.
\end{abstract}




\section{Introduction}
The intent of this article is to summarize some recent progress on the theory of quandles and racks from an algebraic point of view.
Quandles are in general non-associative structures whose axioms correspond to the algebraic axiomatization of the three Reidemeister moves in knot theory \cite{EN}. Quandles and racks appeared in the literature with many different names.  In 1942 Mituhisa Takasaki \cite{Takasaki} introduced   the notion of kei as an abstraction of the notion of symmetric transformation.
The earliest known work on racks is contained in the 1959 correspondence between John Conway and Gavin Wraith who studied
racks in the context of the conjugation operation in a group.   Around 1982,  Joyce \cite{Joyce} and Matveev \cite{Matveev} introduced independently the notion of a quandle.  
Joyce and Matveev associated to each oriented knot a certain quandle in such a way that one of the main problems in knot theory, the problem of equivalence of knots, is transformed into  the problem of isomorphism of quandles.  Since then, quandles  have been of much interest to topologists as well as algebraists.    A Cohomology theory for racks was introduced in \cite{FRS} and was extended to a quandle cohomology in \cite{CJKLS}.  This quandle cohomology was used to define state-sum invariant for knots and paved the road to prove some properties of knots and knotted surfaces, such as invertibility \cite{CJKLS, CEGS} and the minimal triple point numbers of knotted surfaces \cite{SS}.  Quandles were also investigated from an algebraic point of view and  relations to other algebraic structures such as Lie algebras \cite{CCES1, CCES2},  Frobenius algebras and Yang-Baxter equation \cite{CCEKS}, Hopf algebras \cite{And-Grana, CCES2}, quasigroups and Moufang loops \cite{Elhamdadi}, representation theory \cite{EM} and ring theory \cite{BPS}.  For more details on racks, quandles and some other related structures (birack, biquandles) we refer the reader to \cite{EN}.

The outline of this article is as follows.  In Section~\ref{sec2}, we review the basics of quandles and give examples.  Section~\ref{sec3} deals with some notions of representation theory of quandles.  Precisely,  we review the concept of finitely stable quandles, show the existence of non-trivial stabilizing families for the core quandle $Core(G)$ and discuss strong irreducible representations.  In section~\ref{sec4} we relate quandles to ring theory by focusing on the quandle ring of quandles.  We present the key elements 
to show that non trivial quandles are not power associative and we survey doubly transitive quandles and discuss the problem of isomorphism of quandle rings. 

\section{Background on Racks and Quandles}\label{sec2}
Since this paper is meant to be a short overview, we will not discuss here the use of quandles in low dimensional topology and knot theory.   We refer the interested reader, for example, to the book \cite{EN}. 

\noindent
Let $(X, \rRack)$ be a set with a binary operation.  If the set $X$ is closed under the operation $\rRack$, that is for all $x,y$ in $X$, $x\rRack y \in X$, then  $(X, \rRack)$ is called a \emph{magma}.
For $x \in X$, the right multiplication by $x$ is the map $R_x:X \rightarrow X$ given by $R_x(u)=u \rRack x, \forall u \in X$.  Now we give the following definition of a {\em rack}
\begin{df} \cite{EN}
A rack is a magma $(X, \rRack)$ such that  for all $x$ in $X$, the right multiplication $R_x$ is an automorphism of $(X, \rRack)$.
\end{df}
If in addition every element is idempotent, then $(X,\rRack)$ is called a {\em quandle}. 
The following are some examples of racks and quandles
\begin{ex}

\begin{enumerate}
\item
A rack $X$ is {\em trivial} if $\forall x\in X$, $R_x$ is the identity map. 
\item
For any abelian group $G$,  the operation $x \rRack y=2y-x$ defines a quandle structure on $G$ called {\em Takasaki} quandle.  In particular, if $G=  \mathbb{Z}_n$ (integers modulo $n$), 
it is called {\em dihedral quandle} and denoted by $R_n$.

\item

Let $X$ be a module over the ring $\Lambda=\mathbb{Z}_n[t^{\pm 1},s]/(s^2-(1-t)s)$. Then $X$ is a rack with operation $x \rRack y =tx+sy$.  If $t+s \neq 1$, then this rack is not a quandle.  But if $s=1-t$, then this rack becomes a quandle called {\em Alexander} quandle (also called affine quandle).  

\item

The conjugation $x \rRack y=yxy^{-1}$ in a group $G$ makes it into a quandle, denoted $Conj(G)$.

\item

The operation $x \rRack y=yx^{-1}y$ in a group $G$ makes it into a quandle, denoted $Core(G)$.

\item
Let $G$ be a group and let $\psi$ be an automorphism of $G$, then the operation $x \rRack y=\psi(xy^{-1})y$ defines a quandle structure on $G$.  Furthermore, if $H$ is a subgroup of $G$ such that $\psi(h)=h,$ for all $ h \in H$, then the operation $Hx  \rRack Hy=  H \psi(xy^{-1})y$ gives a quandle structure on $G/H$.  
\end{enumerate}

For each $x \in X$, the left multiplication by $x$ is the map denoted by $L_x:X \rightarrow X$ and given by $L_x(y):=x \rRack y$.   A function $f: (X, \rRack) \rightarrow  (Y,\Rack)$ is a quandle {\em homomorphism} if for all $x,y \in X, f(x \rRack y)=f(x) \Rack f(y)$.  Given a quandle $(X,\rRack)$, we will denote by {\rm Aut(X)} the automorphism group of $X$.   The subgroup of {\rm Aut(X)}, generated by the automorphisms $R_x$, is called the {\em inner} automorphism group of $X$ and denoted by {\rm Inn}$(X)$.  The subgroup of {\rm Aut(X)}, generated by  $R_xR_y^{-1},$ for all $ x, y \in X$, is called the $transvection$ group of $X$ denoted by $Transv(X)$.  It is well known \cite{Joyce} that  the $transvection$ group is a normal subgroup of the inner group and the latter group is a normal subgroup of the   automorphism group of $X$.  The quotient group ${\rm Inn(X)}/ {\rm Transv}(X)$ is a cyclic group.  To every quandle, there is a group associated to it which has a {\em universal property} and called {\em enveloping} group of the quandle.  The following is the definition.

\begin{df}
Let $(X, \rRack) $ be a quandle.  The {\em enveloping group} of $X$ is defined by $G_X = F(X) / <x *y = yxy^{-1} , x, y \in X>$, where $F(X)$ denotes the free group generated by $X$.
\end{df}

 Note that there is a natural map $\iota: X \rTo G_X$.   In the following, we show a functoriality property of this group.
\begin{pro}\cite{FR}
Let $(X,\rRack)$ be a quandle and $G$ be any group. Given any quandle homomorphism  $\varphi:X \rTo Conj(G)$. Then, there exists a unique group homomorphism $\widetilde{\varphi}:G_X\to G$ which makes the following diagram commutative
\[
\xymatrix{
X \ar[rr]^{\iota} \ar[d]_{\varphi} &&G_X \ar[d]^{\widetilde{\varphi}}  \\
Conj(G) \ar[rr]_{id} && G
}
\]
\end{pro}
Thus one obtains the natural identification expressing that the functor enveloping from the category of quandles to the category of groups is a {\em left adjoint} to the functor conjugation:  $Hom_{qdle}(X,Conj(G)) \cong Hom_{grp}(G_X,G)$.

The following are some properties and definitions  of some quandles .

\begin{itemize}

\item
A quandle $X$ is {\it involutive}, or a {\it kei}, 
if,  $\forall x \in X,$ ${ R}_x$ is an involution.

\item
A quandle is {\it connected} if ${\rm Inn}(X)$ acts transitively on $X$.

\item

A quandle is {\em faithful} if $x \mapsto {R}_x$ is an injective mapping from $X$ to ${\rm Inn}(X)$.

\item
A {\em Latin quandle} is a quandle such that for each $x \in X$, the left multiplication ${ L}_x$ by $x$ is a bijection. 
That is, the multiplication table of the quandle is a Latin square. 

\item
A quandle $X$ is {\it medial} if 
$(x \rRack y)\rRack (z \rRack w)=(x \rRack z)\rRack (y \rRack w)$  for all  $x,y,z,w \in X$. 
It is well known that a quandle is medial if and only if its tranvection group is abelian.  For this reason, sometimes medial quandles are called {\em abelian}. 
For example, every Alexander quandle is medial.

\item

A quandle $X$ is called $simple$ if the only surjective quandle homomorphisms
on $X$ have trivial image or are bijective.

\end{itemize}
\end{ex}

\section{Elements of Representation Theory of Quandles}\label{sec3}

Recently,  elements of representation theory of racks and quandles were introduced in \cite{EM}.  A certain class of quandles called {\em finitely stable quandles} was introduced as a generalization of the notion of a {\em stabilizer} defined in~\cite{Elhamdadi-Moutuou:Foundations_Topological_Racks}.   A {\em stabilizer} in a quandle $X$ is an element $x\in X$ such that $R_x$ is the identity map.   From the definition of the quandle operation of $Conj(G)$, one sees immediately that an element of $G$ is a stabilizer if and only if it belongs to the center $Z(G)$ of $G$.  On the other hand, an element $x\in Core(G)$ is a stabilizer if and only if $\forall y \in G, (xy^{-1})^2=1$, and thus the  {\em core quandle} $Core(G)$ of a group $G$ has no stabilizers if $G$ has no non-trivial $2$--torsions. This observation suggested that the property of having stabilizers is too strong to capture the identity and center of a group in the category of quandles. This property was then weakened to give the notion of {\em stabilizing families}.%

%

\subsection{Finitely Stable Quandles}\label{fsr}

 
  First we recall the following definition.
\begin{df}\cite{EM}
Let $X$ be a rack or quandle.
 A \emph{stabilizing family of order} $n$  (or $n$--{\em stabilzer})  for $X$ is a finite subset $\{x_1, \ldots, x_n\}$ of $X$ such that the product $R_{x_n}R_{x_{n-1}}\ldots R_{x_1}$ is the identity element in the group $Inn(X)$.
 \end{df}
   
For any positive integer $n$, we define the \emph{$n$--center} of $X$, denoted by $\cal S^n(X)$, to be the collection of all stabilizing families of order $n$ for $X$.  The collection $\cal S(X):= \bigcup_{n\in \bb N}\cal S^n(X)$ of all stabilizing families for $X$ is called the \emph{center} of the quandle $X$. 
%
Now we state the definition of finitely stable quandle.
\begin{df}\label{Bla}\cite{EM}
A quandle $X$ is said to be \emph{finitely stable} if $\cal S(X)$ is non-empty. It will be called $n$--\emph{stable} if it has a stabilizing family of order $n$.
\end{df}

\begin{rmk}
Note that if every element of $X$ is an $n$--stabilzer, then we recover the definition of an $n$--{\em quandle} of Joyce~\cite{Joyce}, Definition 1.5, page 39.
\end{rmk}

\begin{ex}
From Definition~\ref{Bla}, one sees that any finite quandle $X$ is finitely stable.
\end{ex}


\begin{ex}
Consider the integers $\Z$ with the quandle structure $$x\rRack y = 2y - x, \ \ x, y\in \Z.$$ 
Given any family $\{x_1, \ldots, x_n\}$ of integers we get the formula
\[
y \rRack(x_i)_{i=1}^n = 2 \sum_{i=0}^{n-1} (-1)^ix_{n-i} +(-1)^n y, \ \forall y \in \Z.
\] 
One then obtain the stabilizing family $\{x_1,x_1, \cdots, x_{n}, x_{n}\}$ of order $2n$ and thus  $\Z$ admits an infinitely many stabilizing families of even orders.
\end{ex}

Now we show how to obtain non-trivial stabilizing families of the core quandle $Core(G)$.


\begin{pro}\label{pro:Core-2k}\cite{EM}
Let $G$ be a non-trivial group. Then for all even natural number $2k\leq |G|$, there exists a non-trivial stabilizing family of order $2k$ for the \emph{core} quandle $Core(G)$. In particular, if $G$ is infinite, $Core(G)$ admits infinitely many stabilizing families.
\end{pro}

For the stabilizing families of {\em odd order} of the core quandle, we have the following characterization. 

\begin{thm}\label{thm:Core-odd}\cite{EM}
	Let $G$ be a group.  The core quandle of $G$ is $(2k+1)$--stable if and only if all elements of $G$ are $2$--torsions; in other words, $G$ is isomorphic to $\bigoplus_{i\in I}\bb Z_2$, for a certain finite or infinite set $I$.
\end{thm}

Sometimes, a quandle may have no stabilizing family at all as shown in the following example.

\begin{ex}
Let $V$ be a vector space equipped with the quandle structure $x\rRack y = \frac{x+y}{2}.$ Then $\cal S(V)=\emptyset$.
\end{ex}

Let $G$ be a group and $n$ a positive integer. The $n$--\emph{core} of $G$ was defined by Joyce~\cite{Joyce} as the subset of the Cartesian product $G^n$ consisting of all tuples $(x_1, \ldots, x_n)$ such that $x_1\cdots x_n = 1$. Moreover, the $n$-core of $G$ has a natural quandle structure defined by the formula 
\begin{eqnarray}\label{eq1:n-core}
(x_1, \ldots, x_n)\rRack (y_1, \ldots, y_n) := (y_n^{-1}x_ny_1, y_1^{-1}x_1y_2, \ldots, y_{n-1}^{-1}x_{n-1}y_n)
\end{eqnarray}

A notion of the $n$-\emph{pivot} of a group as a generalization of the $n$-core is given by the following definition.
\begin{df}
Let $G$ be a group and let $Z(G)$ be its center. For $n\in \bb N$, we define the \emph{$n$-pivot} of $G$ to be the subset of $G^n$ given by
\[
\cal P^n(G) = \{(x_1, \ldots, x_n)\in G^n \mid x_1\cdots x_n \in Z(G)\}.
\]
\end{df}

It is straightforward to see that the formula~\eqref{eq1:n-core} defines a quandle structure on $\cal P^n(G)$.  Moreover, we have the following proposition.

\begin{pro} \label{pro:core_n} {\rm (\cite{EM})}
Let $G$ be a non-trivial group. Then for all $n\in \bb N$, there is a bijection $$\cal S^n(Conj(G)) \cong \cal P^n(G).$$
In particular, $\cal S^n(Conj(G))$ is naturally equipped with a quandle structure. Furthermore, if $G$ has a trivial center, then $\cal S^2(Cong(G)) = Core(G)$.
\end{pro}

\begin{ex}
If $G$ has trivial center, then for all  $n\geq 2$, $\cal S^n(Conj(G))$ coincides with the $n$--core of $G$. In particular $\cal S^2(Cong(G)) = Core(G)$.
\end{ex}



\subsection{Quandle Actions and Representations} \label{rackaction}

Recall that in~\cite{Elhamdadi-Moutuou:Foundations_Topological_Racks} a \emph{rack action} of a rack $X$ on a space $M$ consists of a map $M\times X \ni (m,x) \mto m\cdot x \in M$ such that
\begin{itemize}
\item[(I)] $\forall x\in X$, the map $ m \mto m\cdot x $ is a bijection of $M$; and 
\item[(II)] $\forall m\in M, x,y\in X,$ we have 
\begin{equation}\label{eq1:rack_action}
(m\cdot x)\cdot y = (m\cdot y)\cdot (x\rRack y).
\end{equation}
\end{itemize}
If $\{x_i\}_{i=1}^s$ is a family of elements in $X$, we will write $m\cdot(x_i)_i$ for 
\[
(\cdots (m\cdot x_1)\cdot \ \cdots )\cdot x_s.
\]

We will require two additional axioms that generalize, in an appropriate way, the concept of group action. Precisely we give the following.

\begin{df}\label{df:rack_action}(\cite{EM})
An action of the rack $(X,\rRack)$ on the space $M$ consists of a map $M\times X \ni (m, x)\mto m\cdot x\in M$ satisfying equation~\eqref{eq1:rack_action} and such that
for all $\{u_1, \ldots, u_s\} \in \cal S(X)$, 
\begin{equation}
  m\cdot (u_i)_i = m\cdot (u_{\s(i)})_i, \ \ \forall m\in M,
\end{equation}
for all cycle $\s$ in the subgroup of $\fr S_n$ generated by the cycle $(2\  \ldots n \ 1).$

\end{df}

\begin{ex}
Let $G$ be a group acting (on the right) on a space $M$.  We then get a rack action of $Conj(G)$ on $M$ by setting $m \cdot g:=mg^{-1}$.  It is easy to check that for all $\{g_i\}_i \in \cal P^n(G)$, we have 
\[
m\cdot (g_i)_i= m \cdot (g_{\sigma(i)})_i,
\]
for all $m \in M$ and $\sigma$ in the subgroup generated by the cycle $(2\  \ldots n \ 1).$

\end{ex}

Notice however that a rack action of $Conj(G)$ does not necessarily define a group action of $G$.

\begin{df}
Let $X$ be a rack acting on a non-empty set $M$ and let $m\in M$. For $\{x_i\}_{i=1}^s\in X$, we define: 
\begin{itemize}
\item[(i)] 
the \emph{orbit} of $\{x_i\}_i$ as 
$
M\cdot (x_i)_i = \{m\cdot (x_i)_{i=1}^s \mid m\in M\}, \text{and}
$

\item[(ii)] the \emph{stabilizer} of $m$ as
$X[m] = \{x\in X \mid m\cdot x = m\}$.
\end{itemize}
\end{df}

Then, it follows immediately from equation~(\ref{eq1:rack_action}), that all stabilizers are subquandles of $X$.

\begin{df}
A quandle action of $X$ on a nonempty set $M$ is {\em faithful} if for any $m \in M$, the map $\varphi_m: X \rightarrow M$ sending $x$ to $m\cdot x$ is injective, i.e. $m \cdot x \neq m \cdot y$ for $x \neq y$.

\end{df}

\begin{lem}\label{lem:faithful_action} {\rm (\cite{EM})}
Let $X$ be a nonempty set. A quandle structure $\rRack$ on $X$ is trivial if and only if $(X,\rRack)$ acts faithfully on a nonempty set $M$ and the action satisfies $(m\cdot x)\cdot y = (m\cdot y)\cdot x$ for all $x,y\in X, m\in M$.
\end{lem}


Now we give a definition of various properties of the action (see \cite{EM}).  
\begin{df}
Let $X$ be a quandle acting on a non-empty set $M$. 
\begin{itemize}
\item[(i)] An {\em approximate unit} for the quandle action is a finite subset $\{t_i\}_{i=1}^r\subset X$ such that $m\cdot (t_i)_i=m$ for all $m\in M$.
\item[(ii)] An element $t\in X$ is an $r$--{\em unit} for the quandle action if the family $\{t\}_{i=1}^r$ is an approximate unit for the rack action.
\item[(iii)] The quandle action of $X$ on $M$ is called {\em $r$--periodic} if each $t\in X$ is an $r$--unit for the quandle action.
\item[(iv)] A quandle action is said to be {\em strong} if every stabilizing family of the quandle is an approximate unit.
\end{itemize}

\end{df}


\begin{ex}
Let $(X, \rRack)$ be a quandle. Then we get a rack action of $X$ on its underlying set by defining $m\cdot x := m\rRack x, \ m, x\in X$. It is immediate that a finite subset $\{t_i\}_i\subset X$ is an approximate unit for this quandle action if and only if it is a stabilizing family for the quandle structure, hence it is a strong action. Furthermore, this action is $n$--periodic if and only if $X$ is an $n$--quandle in the sense of Joyce~\cite{Joyce}, Definition 1.5, p 39.
\end{ex}

Now we define the notion of representation of a quandle (see \cite{EM}). 

\begin{df}

A \emph{representation} of a quandle $X$ consists of a vector space $V$ and a quandle homomorphism 
$
\rho: X\rTo Conj(GL(V)).
$
In other words, we have $\forall x, y\in X, \;\; \rho_{x\rRack y} = \rho_y\rho_x\rho_y^{-1},$ where the isomorphim $\rho(x)$ of the vector space $V$ is denoted by $\rho_x$.
\end{df}

\begin{ex}
Any representation $(V,\rho)$ of a group $G$ induces, in a natural way, a representation of the quandle $Conj(G)$. 
\end{ex}

Let us give an example of quandle representation analogous to the regular representation in group theory. 

\begin{ex}\label{ex:regular}
Let $X$ be a finite quandle and let $\mathbb{C}[X]$ be the vector space of complex valued functions on $X$, seen as the space of formal sums $f = \sum_{x\in X} a_xx$, where $a_x\in \bb C, x\in X$.
We now construct the {\em regular representation} of a quandle $X$ 
\[
\rho : X\rTo Conj(GL(\bb C[X])), \quad \text{given by} \quad \rho_x(f)(y):=f(R_x^{-1}(y)).
\] 
\end{ex}

\begin{df}(\cite{EM})
Let $V$ and $W$ be representations of the quandle $X$. A linear map $\phi: V\rTo W$ is called an {\em intertwining map} of representations
if for each $x\in X$ the following diagram commutes
\[
\xymatrix{
V \ar[r]^{\rho^V_x} \ar[d]_{\phi} & V \ar[d]^{\phi} \\
W \ar[r]^{\rho^W_x} & W.
}
\]
\end{df}
When $\phi$ is an isomorphism then it is called an {\em equivalence of quandle representations}.  If $V$ is a representation of the quandle $X$, a subspace $W\subset V$ is called a {\em subrepresentation} if for any finite family $\{x_i\}_i\subseteq X$, $W\cdot (x_i)_i \subseteq W$.


\begin{df}
A quandle representation $V$ of $X$ is said to be {\em irreducible} if it has no proper subrepresentations.
\end{df}

\subsection{Strong Representations}\label{strong_rep}

\begin{df}\cite{EM}
A {\em strong representation} of $X$ is a representation $V$ such that the quandle action is strong.

We denote by $Rep_s(X)$ the set of equivalence classes of strong {\em irreducible finite dimensional} representations of $X$. 
\end{df}

One can check that as in the group case, $Rep_s(X)$ is an abelian group under tensor product of strong irreducible representations. 

\begin{rmk}
In view of Proposition~\ref{pro:core_n}, we notice that if $G$ is an abelian group, the only strong representation of $Conj(G)$ is the trivial one.
\end{rmk}

\begin{ex}
The regular representation $ (\bb{C}[X], \rho)$ of a rack $X$, defined in Example \ref{ex:regular}, is clearly strong.
\end{ex}

\begin{ex}
Let $(\bb{Z}_3=\{1,2,3\}, \rRack)$, with $x \rRack y=2y-x$, be the dihedral quandle. Define $\rho: \bb{Z}_3 \rTo Conj(GL(\bb{C}^3))$  as the rack representation induced by the reflections on $\bb C^3=Span\{e_1,e_2,e_3\}$
\[
\rho_1= (2 \; 3), \; \rho_2=(1 \; 3),\; \text{and}\; \rho_3=(1 \; 2).
\]
Then $\rho$ is a strong representation of $\bb{Z}_3$ as a quandle, although it is clearly not a group representation. Note, however, that this is a reducible representation since we have the following complete decomposition into irreducible subspaces $\bb C [X]= \bb C<e_1+e_2+e_3> \oplus \; \bb C <e_2-e_1, e_3-e_1>$.
\end{ex}

\begin{ex}\label{ex:involution_rep}
If $X$ is an involutive quandle then, every $x \in X$ is a $2$--stabilizer for $X$. Then every pair $(V,\tau)$, where $V$ is a vector space and $\tau:V\rTo V$ is a linear involution, gives rise to a strong representation $\tilde{\tau}:X\rTo Conj(GL(V))$ by setting $\tilde{\tau}_x(v)=\tau(v)$ for all $x\in X, v\in V$.
\end{ex}


\begin{pro}\cite{EM}
Suppose the quandle $X$ is finite, involutive, and connected. Then, the regular representation of $X$ corresponds to a conjugacy class of the symmetric group $\frak S_n$ where $n={\#X}$. More generally, if $X$ has $k$ connected components, then the regular representation corresponds to $k$ conjugacy classes in $\frak S_n$.  
\end{pro}

%

\begin{thm}\label{thm:irr-stable}
Let $X$ be a finite quandle. Then every irreducible strong representation of $X$ is either trivial or finite dimensional.
\end{thm}

We shall observe that in the case of involutive quandles the above result becomes more precise with regard to the dimension of the irreducible representations.   We then have the following theorem (see \cite{EM})

\begin{thm}\label{thm:strong_rep_connected}
Every irreducible strong representation of an involutive connected finite quandle $X$ is one-dimensional.
\end{thm}

\section{Quandles and Ring Theory}\label{sec4}

 In \cite{BPS}, a theory for quandle rings was proposed for quandles by analogy to the theory of group rings for groups.  In \cite{EFT}, quandles were investigated from the point of view of ring theory.  Precisely, some basic properties of quandle rings were investigated and some open questions which were raised in \cite{BPS} were solved.  In particular, examples of quandles were provided for which the quandle rings $\mathbf{k}[X]$ and $\mathbf{k}[Y]$ are isomorphic, but the quandles $X$ and $Y$ are not isomorphic.  

\subsection{Power Associativity of Quandle Rings}\label{Powerassoc}

 We associate to every quandle $(X, \rRack)$ and an associative ring $\mathbf{k}$ with unity, a \textit{nonassociative} ring $\mathbf{k}[X]$.  First, we start by recalling the following Definition and Proposition. 
 
 \begin{dfpro}
 Let $(X, \rRack)$ be a quandle and $\mathbf{k}$ be an associative ring  with unity.
 Let $\mathbf{k}[X]$ be the set of elements that are uniquely expressible in the form $\sum_{x \in X }  a_x x$, where $x \in X$ and $a_x=0$ for almost all $x$.  Then the  set $\mathbf{k}[X]$ becomes a ring with the natural addition and the multiplication given by the following, where $x, y \in X$ and $a_x, a_y \in \mathbf{k}$,
\[   ( \sum_{x \in X }  a_x x) \cdot ( \sum_{ y \in X }  b_y y )
=   \sum_{x, y \in X } a_x b_y (x \rRack y) . \]
\end{dfpro}

 In \cite{BPS}, power associativity of dihedral quandles was investigated and the question of determining the conditions under which the  quandle ring $R[X]$ is power associative was raised.  In this section we give a complete solution to this question.  Precisely, we prove that quandle rings are never power associative when the quandle is non-trivial and the ring  has characteristic zero.	
  
Let first let's recall the following definition from \cite{Albert} 
\begin{df}
	A ring $\mathbf{k}$ in which every element generates an associative subring is called a \textit{power-associative} ring. 
\end{df}

\begin{ex}{\rm
		Any alternative algebra  is power associative.  Recall that an algebra $A$ is called \textit{alternative} if $x \cdot (x \cdot y)= (x \cdot x )\cdot y$ and  $x \cdot (y \cdot y)= (x \cdot y )\cdot y, \forall x, y \in A$,  (for more details see \cite{EMakh}).  
	} \end{ex}
	In non-associative algebras,  the product $x \cdot x$ is also denoted by  $x^2$, as in the associative case.

		It is well known \cite{Albert} that a ring $\mathbf{k}$ of characteristic zero is power-associative if and only if, for all $ x \in \mathbf{k},$
		\begin{eqnarray}
		x^2\cdot x &=& x\cdot x^2, \label{three} 
		\end{eqnarray}
		and
		\begin{eqnarray}
		x^2 \cdot x^2 &=& (x^2 \cdot x) \cdot x. \label{four}
		\end{eqnarray}
\begin{thm} \cite{EFT}
Let $\mathbf{k}$ be a ring  of characteristic zero and let $(X,\rRack)$ be a non-trivial quandle.  Then the quandle ring $\mathbf{k}[X]$ is not power associative.
\end{thm}
The proof of this theorem is based on the use of the two identities~(\ref{three}) and~(\ref{four}).  The proof consists of two steps: 

\begin{enumerate}
\item 
 show that if there exist $x, y\in X$ such that $x\neq y$ and $x\rRack y = y\rRack x$, then $\mathbf{k}[X]$ is not power associative.
\item 
Assume  $ \forall x, y\in X$ such that $x\neq y$, we have  $x\rRack y \neq y\rRack x$. If $\mathbf{k}[X]$ is power associative, then $x\rRack y= x$.
\end{enumerate}
A complete proof can be found in \cite{EFT}.

\subsection{Higher Transitivity in Groups and Quandles }
First we recall the notion of  $n$-transitivity of groups \cite{DM} and use it in the context of quandles.  Let $G$ be a group acting (the right) on a set $X$, thus $G$ acts naturally on any cartesian product $X^n$ of $X$ by $(x_1, \cdots, x_n)g= (x_1 g, \cdots, x_n g) $.  Note that the subset $X_n:=\{(x_1, \cdots, x_n) \in X^n | \; \forall i \neq j, \;x_i \neq x_j  \}$, of $X^n$, is $G$-invariant.  We then state the following definition of $n$-transitive action.
\begin{df}\cite{DM}
Let $G$ be a group, $X$ be a non-empty set and $n$ be a positive integer.  The group $G$ acts $n$-transitively on $X$ if it acts transitively on the subset $X_n$ of $X^n$.
\end{df}
It is clear from this definition that an $n$-transitive action implies an $(n-1)$-transitive action.  It is also well known that $Sym(X)$ acts $n$-transitively on $X$ when $n \leq \#X$.

Now we state the following definition of $n$-transitive quandles.

\begin{df}\cite{McCarron1}
A finite quandle $X$, with cardinality $|X|>2$, is called $n$-transitive if the action of the inner group $Inn(X)$ is $n$-transitive action on $X$.
\end{df}
The following is the definition of {\em primitive} action.
\begin{df}\cite{DM}
Let $G$ be a group acting on (the right) on a finite non-empty set $X$.  We say that $G$ preserve an equivalence relation $\sim$ on $X$ if the following condition is satisfied: 
\[
x \sim y \;  \text{if and only if} \;xg \sim yg.
\]
The action of $G$ on $X$ is said to be primitive if the group $G$ acts transitively on $X$ and if $G$ preserves no non-trivial partition of $X$.
\end{df}

It is well known \cite{DM} that any $2$-transitive group is primitive, but the converse is not true in general.  Corollary 1.5A on page 14 of \cite{DM} states that a group $G$ is primitive if and only if each point stabilizer is a maximal subgroup of $G$.  Thus, the study of finite primitive permutation groups is equivalent to the study of the maximal subgroups of finite groups.

McCarron showed in Proposition 5 of \cite{McCarron1} that if $n \geq 2$ and $X$ is a finite
$n$-transitive quandle with at least four elements, then $n= 2$.  Furthermore, he showed that the dihedral quandle of three elements is the only $3$-transitive quandle, and there are no $4$-transitive quandles with at least four elements.\\
Now we turn to the notion of quandles of {\em cyclic type} and we will give their relastionship to $2$-transitive quandles (called also {\em two point-homogeneous} in \cite{Tamaru}).

\begin{df}
A finite quandle $X$ is of {\em cyclic type} if for each $x \in X,$ the permutation $R_x$ acts on $X \setminus \{x\}$ as a cycle of length $|X| - 1$, where
$|X|>2$ denotes the cardinality of X. 
\end{df}

\begin{ex}
The following list of quandles of cyclic type is taken from table 1 in \cite{KTW}.
\begin{itemize}
\item
Among the quandles of order three, the dihedral quandle is the only one of cyclic type.

\item

Among the quandles of order four, the quandle $X=\{1,2,3,4\}$ with $R_1=(234)$, $R_2=(143)$, $R_3=(124)$ and $R_4=(132)$ is the only one of cyclic type.

\item

On the set of five elements  $\{1,2,3, 4, 5\}$, there are exactly two quandles $X$ and $Y$ of cyclic type:\\
 $X$ with $R_1=(2345)$, $R_2=(1354)$, $R_3=(1425)$ $R_4=(1532)$ and $R_5=(1243)$, and \\
  $Y$ with $R_1=(2345)$, $R_2=(1435)$, $R_3=(1542)$ $R_4=(1253)$ and $R_5=(1324)$.\\
\end{itemize}
\end{ex}
The classification of quandles of {\em cyclic type} was investigated in \cite{KTW}.   Their main theorem states that the isomorphism classes of cyclic quandles of order $n$ are in one-to-one correspondence with the permutations of the symmetric group $\fr S_n$ satisfying certain conditions.  Furthermore, they obtained the classification of the quandles of cyclic type with cardinality up to order $12$. 

They also stated the following conjecture: "A quandle with at lease $3$ elements is of cyclic type if and only if its cardinality is a power of a prime".  This conjecture was later proved in \cite{Vend} (see Theorem $1$ and also Corollary $2$).  These types of quandles appeared earlier in \cite{LR}, where they were called quandles with {\em constant profile} $(\{1, n-1\}, \ldots,  \{1, n-1\})$.   

In \cite{Tamaru}, the author defined a quandle $X$ to be two-point homogeneous if the inner automorphism group $Inn(X)$ acts doubly-transitive on $X$.   Furthermore, the author classified two-point homogeneous quandles of prime cardinality by proving that they are linear Alexander quandles.  This was generalized in \cite{Wada, Vend} to the case of prime powers giving the following result: A finite quandle $(X,\rRack)$ is $2$-transitive \emph{if and only if} it is of cyclic type,  \emph{if and only if}  it is isomorphic to an Alexander quandle $(\mathbb{F}_q, \rRack)$ over the finite field $\mathbb{F}_q $ with operation $x \rRack y= ax+(1-a)y$ where $a$ is a primitive element (i.e. a generator of the multiplicative group $\mathbb{F}_q ^{\times}$) .

\par
The notion of $2$-transitive quandles was also used recently in \cite{EFT} to study quandles from a ring theoretical approach.



\subsection{Isomorphisms of Quandle Rings}

 In \cite{EFT}, a notion of partition-type of quandles was introduced and it was shown that if the quandle rings $\mathbf{k}[X]$ and $\mathbf{k}[Y]$ are isomorphic and the quandles $X$ and $Y$ are orbit $2$-transitive, then $X$ and $Y$ are of the same partition type.  Now let us recall the definition of partition type.
	
\begin{df}(\cite{EFT})
Let $X$ be a finite quandle of cardinality $n$. The \textit{partition type} of $X$ is  $\lambda = (\lambda_1, \hdots, \lambda_n)$ with  $\lambda_j$ being the number of orbits of cardinality $j$ in $X$. 
\end{df}	

\begin{ex}
Let $X=\{1,2,3,3,5,6,7\}$ be a quandle with the orbit decomposition $X=\{1\} \amalg \{2,3\} \amalg \{ 4,5,6\} \amalg \{7\}$.  Then the  partition type of the quandle $X$ is $\lambda=(2,1,1,0,\hdots)$.
\end{ex}
In order to state the next result about the same partition type for quandles with isomorphic quandle rings, we need the following theorem which gives the list of subgroups $G\leqslant  S_n$ for which the representation $V_{st}$ is irreducible, where $V_{st}\subset V= \mathbf{k}[X]$ be the subspace orthogonal to the vector $v_{triv}=\sum\limits_{x\in X_i} x$. 
The groups $(i)-(v)$ are respectively the affine and projective general semilinear groups, projective semilinear unitary groups, Suzuki and Ree groups. We refer the reader to \cite[Chapter 7]{DM} for definitions of these groups. 
\par
\begin{thm} \cite{EFT}
\label{TransClassif}
\begin{enumerate}
\item[(a)]

 If $char(\mathbf{k})=0$  then $V_{st}$ is an irreducible representation of the subgroup $G < S_n$ if and only if  $G$ is $2$-transitive and $A_n\not\leq G$. 
 \item[(b)] 
 In case $char(\mathbf{k})=p>3$  then $V_{st}$ is an irreducible representation of the subgroup $G < S_n$ if and only if  $G$ is $2$-transitive and $A_n\not\leq G$ except 
 \begin{enumerate}
 \item[(i)] $G \leqslant  A\Gamma L(m, q)$, and $p$ divides $q$;
 \item[(ii)] $G \leqslant  P\Gamma L(m, q)$, $m\geq 3$ and $p$ divides $q$;
 \item[(iii)] $G \leqslant  P\Gamma U(3, q)$, and $p$ divides $q+1$;
 \item[(iv)] $G \leqslant  Sz(q)$, and $p$ divides $q+1+m$, where $m^2=2q$;
 \item[(v)] $G \leqslant  Re(q)$, and $p$ divides $(q+1)(q+1+m)$, where $m^2=3q$.
\end{enumerate} 
\end{enumerate}
\end{thm}

The following theorem was proved in \cite{EFT}.
\begin{thm}
\label{ISO}
Assume $char(\mathbf{k})\neq 2,3$. If the quandle rings $\mathbf{k}[X]$ and $\mathbf{k}[Y]$ are isomorphic and the quandles $X$ and $Y$ are orbit $2$-transitive  ($G_{X_i}$'s are not among the groups from $(i)-(v)$ in Theorem~\ref{TransClassif} in case $char(\mathbf{k})=p>3$), then $X$ and $Y$ are of the same partition type.
\end{thm}

The following two examples were provided in \cite{EFT} as answers to \textit{Question}  $7.4$ of \cite{BPS} on the existence of two nonisomorphic quandles X and Y with isomorphic quandle rings.   
 \begin{ex}
\label{CounterEx1}
Let $\mathbf{k}$ be a field with $char(\mathbf{k})=3$ and let $X=Y=\{1,2,3,4\}$ as sets.  Define a quandle structure on the set $X$ by $R_1$ and $R_2$ are the identity maps and $R_3$ and $R_4$ are both equal to the transposition $(1\; 2)$.  Define a quandle operation on the set $Y$ by $R_3$ being equal to the transposition $(1\; 2)$ and all others $R_i$ are the identity map.	  Note that the orbit decomposition $\{1,2\}  \amalg \{3\} \amalg \{4\} $ is the same for both quandles $X$ and $Y$.  An isomorphism $f$ between $X$ and $Y$ must satisfy $f(1) \in \{1,2\}$.  The sets $\{x, 1\rRack x=1\}$ and $\{x, 2\rRack x=2\}$ have cardinality $2$ as subsets of $X$ but cardinality $3$ as subsets of $Y$, thus the quandles $X$ and $Y$ are not isomorphic.
The ring isomorphism  is given by $\varphi:\mathbf{k}[X]\overset{\sim}{\rightarrow}\mathbf{k}[Y]$, where $\varphi=\left(\begin{array}{cccc}
1 & 0 & 0 & 1\\
0 & 1 & 0 & 1\\
0 & 0 & 1 & 1\\
0 & 0 & 0 & 1\\
\end{array}\right)$.
\end{ex}

\begin{ex}
\label{CounterEx2}
Let $\mathbf{k}$ be a field with {\em characteristic zero} and  assume $X=Y=\{1,2,3,4, 5,6,7\}$ as sets.  Define a quandle structure on the set $X$ by $R_5=(1\; 2)$,  $R_6=(1\; 2) (3 \; 4)$ and all the other $R_i$ , $i \neq 5,6,$ are the identity map.  Define a quandle operation on the set $Y$ by $R_5=(1\; 2)$,  $R_6=(3 \; 4)$ and all the other $R_i$ , $i \neq 5,6,$ are the identity map.  A similar argument to the one given in the previous example shows that $X$ and $Y$ are not isomorphic as quandles.  One family of ring isomorphisms  is given by $\varphi:\mathbf{k}[X]\overset{\sim}{\rightarrow}\mathbf{k}[Y]$, where 
\[\varphi=\left(\begin{array}{ccccccc}
1 & 0 & 0 & 0 & 0 & 0 & 0\\
0 & 1 & 0 & 0 & 0 & 0 & 0\\
0 & 0 & 1 & 0 & 0 & 0 & 0\\
0 & 0 & 0 & 1 & 0 & 0 & 0\\
0 & 0 & 0 & 0 & 1 & 1 & 0\\
0 & 0 & 0 & 0 & 0 & 1 & 0\\
0 & 0 & 0 & 0 & 0 & -1 & 1\\
\end{array}\right).\]
\end{ex}
There is an easy way to obtain other examples from the one we just provided (see \cite{EFT} for more details).

\noindent
\textbf {Acknowledgement:}
The author would like to thank Neranga Fernando,  El-ka\"ioum M. Moutuou and Boris Tsvelikhovskiy for many fruitful suggestions which improved the paper.

\end{document}